\documentclass[11pt]{article}
\usepackage{amsmath}
\usepackage{amssymb}
\usepackage{amsthm}
\usepackage{cases}
\usepackage{cite}
\usepackage{mathrsfs}
\usepackage{latexsym}
\usepackage{wasysym}
\usepackage[misc]{ifsym}
\usepackage{indentfirst}
\usepackage{sectsty}
\usepackage{times}

\numberwithin{equation}{section}

\newtheorem{theorem}{Theorem}[section]
\newtheorem{proposition}[theorem]{Proposition}
\newtheorem{definition}[theorem]{Definition}
\newtheorem{corollary}[theorem]{Corollary}
\newtheorem{lemma}[theorem]{Lemma}
\newtheorem{remark}[theorem]{Remark}
\newtheorem{example}[theorem]{Example}

\newcommand{\cA}{{\mathcal A}}

\newcommand{\cC}{{\mathcal C}}
\newcommand{\cD}{{\mathcal D}}

\newcommand{\cI}{{\mathcal I}}
\newcommand{\cK}{{\mathcal K}}
\newcommand{\cO}{{\mathcal O}}

\newcommand{\sE}{{\mathscr E}}

\newcommand{\sM}{{\mathscr M}}

\def\R{\mathbb{R}}

\def\N{\mathbb{N}}
\def\1{\mathbb{1}}
\newcommand{\di}{\mathrm{d}}

\newcommand{\mB}{\mathrm{B}}
\newcommand{\cat}{{\rm cat}}

\def\disp{\displaystyle}

\def\bc{\begin{center}}
\def\ec{\end{center}}
\def\be{\begin{equation}}
\def\ee{\end{equation}}
\def\bea{\begin{eqnarray}}
\def\eea{\end{eqnarray}}
\def\ba{\begin{array}}
\def\ea{\end{array}}
\def\benu{\begin{enumerate}}
\def\eenu{\end{enumerate}}
\def\bt{\begin{theorem}}
\def\et{\end{theorem}}
\def\bl{\begin{lemma}}
\def\el{\end{lemma}}
\def\bco{\begin{corollary}}
\def\eco{\end{corollary}}
\def\bn{\begin{numcases}}
\def\en{\end{numcases}}
\def\br{\begin{remark}}
\def\er{\end{remark}}
\def\bd{\begin{definition}}
\def\ed{\end{definition}}
\def\bp{\begin{proposition}}
\def\ep{\end{proposition}}
\def\bo{\begin{proof}}
\def\eo{\end{proof}}
\def\bx{\begin{example}}
\def\ex{\end{example}}

\def\pa{\partial}
\def\al{\alpha}
\def\De{\Delta} \def\de{\delta}
\def\na{\nabla}

\def\ve{\varepsilon}

\def\vp{\varphi}
\def\w{\omega}\def\W{\Omega}
\def\gam{\gamma}

\def\~{\widetilde}

\def\ol{\overline}
\def\ul{\underline}

\def\Cup{\bigcup}

\def\ra{\rightarrow}

\def\stac{\stackrel}
\def\8{\infty}
\def\X{\times}

\def\mb{\mbox}
\def\di{{\rm d}}

\def\es{\emptyset}

\def\sm{\setminus}

\def\ss{\subset}

\def\Hs{\hspace{0.8cm}}
\def\hs{\hspace{0.4cm}}
\def\Vs{\vskip10pt}
\def\vs{\vskip5pt}

\def\[{\left[}
\def\]{\right]}
\def\({\left(}
\def\){\right)}

\title{On Relative Category and Morse Decompositions\\ for Infinite-Dimensional Dynamical Systems\thanks{This work was supported by NSFC grants 11801190.}}
\author{Jintao Wang\thanks{Corresponding author, E-mail address: wangjt@wzu.edu.cn}\\
{\small Department of Mathematics, Wenzhou University, Wenzhou 325035, China}\\
Desheng Li\thanks{E-mail address:lidsmath@tju.edu.cn}\\
{\small School of Mathematics Tianjin University, Tianjin 300072, China}
}

\date{}
\begin{document}
\maketitle
\normalsize
\begin{center}
\begin{minipage}{14cm}
\noindent{\bf Abstract.} We employ the relative category to develop relations between the Wa\.zewski pair $(N,E)$ and the Morse decomposition of the maximal invariant set in $\ol{N\sm E}$ for infinite-dimensional dynamical systems.
Via these relations, we can detect connecting trajectories between Morse sets and obtain a dynamical-system version of critical point theorem with relative category.
\Vs

\Vs
\end{minipage}
\end{center}
\pagenumbering{arabic}
\section{Introduction}

In the field of dynamical systems, invariant sets are of great significance and particular interests, in that they can determine and describe much of the long-term dynamics of a system.
Equilibria, (almost) periodic solutions, homoclinic (heteroclinic) orbits and attractors are typical examples of compact invariant sets.
It is therefore of great importance to study the existence, number and the location of compact invariant sets for a given dynamical system.

Early in 1940s, Wa\.zewski introduced the famous Wa\.zewski's Retract Theorem (\cite{W47R,W47A}) to give the existence of invariant sets.
Roughly speaking, this theorem states that for a given flow and a closed subset $N$ (a Wa\.zewski set) of the phase space, the existence of a solution entirely contained in $N$ can be deduced from the assumption that the exit set $N^-$ of $N$ is not a deformation retract of $N$.
The Wa\.zewski's Retract Theorem appeared to be a powerful way to develop topological methods to study invariant sets.
The well-known Conley index theory (\cite{Con78}) was originally inspired by this theorem, which was generalized to the infinite-dimensional case by Rybakowski (\cite{R80}) and shape index theory by Robbin and Salamon (\cite{RS88}, see more generally, \cite{S03,S07,WLD15,WLD20}).
Li, Shi and Song in \cite{LSS15} used the Wa\.zewski pair (a generalization of the Wa\.zewski set) to develop a dynamical-system version of the linking theorem and mountain pass theorem.
Among these explorations of relevant topics, Lusternik-Schnirelmann category (L-S category for short) has also been introduced to dynamical systems via (Conley, shape) index pairs (see \cite{P94,R00,S00,W17}) in recent decades.

L-S category, as a numerical topological invariant, was introduced by Lusternik and Schnirelmann \cite{LS34} in the course of research into the calculus of variations.
It is shown to provide important information about the existence of critical points.
Since then, L-S category has received a great deal of study over the years \cite{B78,CLOT03,FLR94,J95} and applied extensively to the existence of critical points in calculus of variations.
In recent decades, L-S category has also been applied to the study of invariant sets of dynamical systems.
Po\'zniak introduced a new category as a modification of L-S category in \cite{P94}, in which a relation of categories was presented between an isolated invariant set and its Conley index pair.
Based on the index pairs, Sanjurjo studied the flow on locally compact metric spaces and applied L-S category in the shape theory (\cite{B78}) to the Morse decompositions of an isolated invariant set to detect the existence of connecting orbits between Morse sets (\cite{S00}).
Razvan in \cite{R00} gave a generalization of Po\'zniak's work \cite{P94} and obtained an relationship of the categories of the Conley index and the Morse decomposition of a compact isolated invariant set for a flow to study a lower number of the equilibria.
Weber studied in \cite{W17} the number of isolated critical points for a downward gradient flow on a (smooth) closed manifold via Conley (index) pairs and L-S category.
Nevertheless, those works (\cite{P94,R00,S00,W17}) above about L-S category in dynamical systems were all in the framework of continuous flows, some of which can be applied to the semiflows with specific conditions (e.g. two-sided on the unstable manifold of a given invariant set).

In this present work, we consider the local semiflows on infinite-dimensional metric spaces and study the relation between the relative category of a Wa\.zewski pair $(N,E)$ and the (L-S) category of the Morse decomposition of the maximal invariant set in $N\sm E$.
We shall generalize the corresponding consequences for flows on locally compact metric spaces stated in \cite{S00}.
The definitions of L-S category and relative category we adopt here are referred to \cite{C86,CLOT03,FLR94,J95,W96}; see Section 4.
We use the notation of relative category given in \cite{FLR94}, i.e., for a closed pair $(N,E)$ and $A\ss N$ such that $E\ss A\ss N$, ${\rm cat}_{N,E}(A)$ denotes the {\em category of $A$ in $N$ relative to $E$}.
The classical {\em L-S category} ${\rm cat}_{N}(A)$ is indeed equal to ${\rm cat}_{N,\es}(A)$.

Now we give a detailed description of our work.
Let $\Phi$ be a local semiflow on a complete metric space $X$.
Given a compact invariant set $K$ with its Morse decomposition $\sM=\{M_1,\cdots,M_n\}$, as long as $K$ is contained in an ANE $N$ (see the definition of ANE in Section 4), we have the following consequence,
$$\cat_N(K)\leqslant\sum_{i=1}^n\cat_N(M_i).$$
This is a generalization of \cite[Lemma 4]{S00}.
Since the topological structure of compact invariant sets is usually complicated, it leads to great difficulty in the calculation of the corresponding categories.
This prompts us to study more about this topic around the invariant set $K$ in the theory of relative category.

Let the closed pair $(N,E)$ be a {\em Wa\.zewski pair}, i.e., $E$ is an exit set of $N$ and moreover, since $X$ can be infinite-dimensional, the set $\ol{N\sm E}$ is imposed an appropriate compact condition --- {\em strong admissibility}.
As to establish the relation between $\cat_{N,E}(N)$ and the invariant sets in $N\sm E$, we also suppose the Wa\.zewski pair $(N,E)$ satisfies the {\em transversality} condition, which ensures that each point lying on the boundary of $E$ in $N$ will not stay on this boundary within an arbitrarily short time under the action of semiflow, namely, Definition \ref{de3.2}.
Our main theorem is the following one.

\bt\label{th1.1}Let $\Phi$ be a local semiflow on $X$ and $(N,E)$ be a transversal Wa\.zewski pair for $\Phi$ with $E\subset N$.
Let $\{M_1,\cdots,M_n\}$ be a Morse decomposition of the maximal invariant set in $\ol{N\sm E}$.
Suppose that $N$ is an ANE.
Then
\be\cat_{N,E}(N)\leqslant\sum_{i=1}^n\cat_N(M_i).\label{1.1}\ee\et

This consequence generalizes those relevant ones of \cite{P94,R00,S00,W17}.
Theorem \ref{th1.1} can also help to establish a corresponding relation between relative category and Conley index; see Remark \ref{re4.13} below.
When the system $\Phi$ is gradient, \eqref{1.1} can be also applied to give a lower-bound estimate for the number of equilibria; see Remark \ref{re4.15}.
This is also an elementary generalization of critical point theorem with relative category (see\cite[Theorem 5.19]{W96}, a sort of minimax theorems) in dynamical systems.

This paper is organized as follows.
In Section 2, we present the preliminaries, including some basic topological concepts, dynamical systems on normal Hausdorff topological spaces, attractors and Morse decompositions.
In Section 3, we introduce the Wa\.zewski pair and the concept --- transversality for a closed pair, and develop some results on the quotient flow.
We state and prove the main theorem in the 4th section.
In Section 5, we give some extensions (or applications) of the categories for Morse decompositions.

\section{Preliminaries}

Let $X$ be a topological space and $A,\,B\ss X$ with $A\ss B$.
We denote by $\ol{A}$ the {\em closure} of $A$ in $X$ and by ${\rm int}_B(A)$ the {\em interior} of $A$ in $B$, i.e., the maximal open subset of $B$ contained in $A$.
A set $U$ is called an ({\em open}, {\em closed}) {\em neighborhood} of $A$ in $B$, if $U$ is (open, closed) in $B$ and there is an open subset $\cO$ of $B$ such that $A\ss\cO\ss U$.

The set $A$ is said to be {\em sequentially compact}, if each sequence $x_n$ in $A$ has a subsequence converging to a point $x\in A$.
It is a basic knowledge that if $X$ is a metric space, then sequential compactness coincides with compactness.

Let $X$ be a Hausdorff topological space.
\bd A {\bf local semiflow} $\Phi$ on $X$ is a continuous map $\Phi:\,\cD(\Phi)\ra X$, where $\cD(\Phi)$  is an open subset of $\R^+\X X$, and $\Phi$ enjoys the following properties:
\vs
(1) for each $x\in X$, there exists $0<T_x\leqslant\8$ such that
$$(t,\,x)\in\cD(\Phi)\Longleftrightarrow 0\leqslant t<T_x;$$

(2) $\Phi(0,\,x)=x$ for all $x\in X$;

(3) if $(t+s,\,x)\in\cD(\Phi)$, where $t$, $s\in \R^+$, then $\Phi(t+s,\,x)=\Phi(t,\,\Phi(s,\,x))$.
\vs
The number $T_x$ in (1) is called the {\bf maximal existence time} of $\Phi(t,x)$.
In the case when $\cD(\Phi)=\R^+\X X$,  we call $\Phi$  a {\bf global semiflow}.
\ed

Let $\Phi$ be a given local semiflow on $X$. For notational convenience, we will  rewrite $\Phi(t,\,x)$ as $\Phi(t)x$. Given a subset $N$ of $X$, we say $\Phi$ {\em does not explode} in $N$, if $T_x=\8$, whenever $\Phi(t)x\in N$ for all $t\in[0,T_x)$.

A subset $N$ of $X$ is said to be {\em admissible}, if for arbitrary sequences $x_n\in N$ and $t_n\ra+\8$ with $\Phi([0,t_n])x_n\subset N$ for all $n$,  the sequence of the end points $\Phi(t_n)x_n$ has a convergent subsequence. The subset $N$ is {\em strongly admissible} if $N$ is admissible and $\Phi$ does not explode in $N$.
The strong admissibility condition can be viewed as the asymptotic compactness of dynamical systems (see \cite{T98}), but locally.

A {\it solution} ({\it trajectory}) on an interval $J\subset \R$ is a map $\gam:\,J\ra X$ satisfying
$$\gam(t)=\Phi(t-s)\gam(s),\Hs\mb{for all }\,s,\,t\in J,\,\,s\leqslant t.$$
A {\em full solution} $\gam$ is a solution defined on the whole line $\R$. If $x\in X$ is such that $\Phi(t)x=x$ for all $t\geqslant0$, we say $x$ is an {\em equilibrium}.

The {\em $\w$-limit set} and {\em $\al$-limit set} of a solution $\gam$ are defined as follows: if $\gam$ is defined on an interval containing $[0,\8)$, it is defined that
$$\w(\gam)=\{y\in X:\,\mb{there exists }t_n\ra\8\mb{ such that }\gam(t_n)\ra y\};$$
if $\gam$ is defined on an interval containing $(-\8,0]$, it is defined that
$$\al(\gam)=\{y\in X:\,\mb{there exists }t_n\ra-\8\mb{ such that }\gam(t_n)\ra y\}.$$
For an $x\in X$ with $T_x=\8$, we define $\w(x)=\w(\gam)$ with $\gam(t)=\Phi(t)x$ for every $t\geqslant0$.

Given an invariant set $K\ss N\ss X$, we define the {\em local stable and unstable manifold}, $W^s_N(K)$ and $W^u_N(K)$ of $K$ in $N$ as follows:
$$W^s_N(K):\,=\Cup_{\w(\gam)\ss K}\{\gam(t):\,\gam([0,\,\8))\ss N,\,t\in[0,\,\8)\},$$
$$\mb{and}\hs W^u_N(K):\,=\Cup_{\al(\gam)\ss K}\{\gam(t):\,\gam((-\8,\,0])\ss N,\,t\in(-\8,\,0]\},$$
where $\gam$ is a solution and $\w(\cdot)$ and $\al(\cdot)$ are limit sets. If $N=X$ is the whole phase space, we simply write $W^s(K)=W^s_X(K)$ and $W^u(K)=W^u_X(K)$.

Let $M$ and $B$ be two subsets of $X$. We say that {\em $M$ attracts $B$}, if $T_x=\8$ for all $x\in B$ and moreover,
for each neighborhood $U$ of $M$ there exists $T>0$ such that $$\Phi(t)B\ss U,\Hs t>T.$$

A nonempty sequentially compact invariant set $\cA\ss X$ is said to be an {\em attractor} of $\Phi$,
if it attracts a neighborhood $U$ of $\ol{\cA}$ and $\cA$ is the maximal sequentially compact invariant set in $U$.
It can referred to \cite[Remark 2.1]{WLD20} about the discussion and comparison of this definition and the previous ones of attractor.

Let $\cA$ be an attractor. Set
 $$\W(\cA)=\{x\in X:\,\cA\mb{ attracts }x\}.$$
$\W(\cA)$ is called the {\em region of attraction} of  $\cA$.
One can easily verify that $\W(\cA)$ is open; moreover, $\cA$ attracts each compact subset of $\W(\cA)$ (see \cite{LWX17}).
In the case when  $\W(\cA)=X$,  we simply  call $\cA$  the {\em global attractor} of $\Phi$.
%
\Vs

For the reader's convenience, we recall briefly the Morse decompositions of invariant sets for dynamical systems on topological systems (see \cite{LWX17} or more classically in \cite{Con78,R80}).

Let $K$ be a compact invariant set. Then the restriction $\Phi|_K$ of $\Phi$ on $K$ is a semiflow on $K$.
A set  $A\ss K$ is called an {\em attractor of $\Phi$ in $K$}, if it is an attractor of $\Phi|_K$.
We denote by $\W_K(A)$ the region of attraction of $A$ in $K$ for $\Phi|_K$.

\bd\label{de2.4} Let $K$ be a compact invariant set. An ordered collection
$\sM=\{M_1,\,\cdots,\,M_n\}$
of subsets $M_k\ss K$ is called a {\bf Morse decomposition} of $K$, if there exists an increasing sequence
$\es=A_0\subsetneq A_1\subsetneq\cdots\subsetneq A_n=K$ of attractors in $K$ such that
$$M_k=A_k\sm\W_K(A_{k-1}),\hs 1\leqslant k\leqslant n.$$

The attractor sequence of $A_k$ ($k=0,\,1,\,\cdots,\,n$) is often called the {\bf Morse filtration} of $K$, and each $M_k$ is called a
{\bf Morse set} of $K$.
\ed

\br If $K$ is an attractor in $X$, each attractor $A$ in $K$ is also an attractor in $X$ for $\Phi$.
Moreover, $\W_K(A)=\W(A)\cap K$.
\er

\section{Wa\.zewski Pairs and Quotient Flows}

In this section, we recall the Wa\.zewski pair and some properties of quotient flows (see \cite{WLD15}), and develop some new related conclusions for this paper. We always assume the phase space $X$ to be a complete metric space with the metric $\di(\cdot,\cdot)$.
\subsection{Wa\.zewski pairs and quotient flows}
Let $A$, $N$ be subsets of $X$. $A$ is said to be {\em $N$-positively invariant},
if $\Phi([0,\,t])x\ss N$ ensures $\Phi([0,\,t])x\ss A$, for every $x\in N\cap A$ and $t\geqslant0$.
When $N=X$, $N$-positive invariance is exactly the positive invariance.
For an arbitrary subset $N\ss X$,  define a function $t_N:\,X\ra\R^+\cup\{\8\}$ as
$$t_N(x)=\inf\{t\geqslant0:\,\mb{either }t\geqslant T_x\mb{ or }\Phi(t)x\not\in N\},\Hs\mb{for all }x\in X.$$
Note that for each $x\in N$, $t_N(x)$ is the supremum of the time $t$ such that $\Phi([0,\,t])x\subset N$.

Let $N$, $E$ be subsets of $X$. We say $E$ is an {\em exit set} of $N$, if
\benu\item[(1)] $E$ is $N$-positively invariant;
\item[(2)] for every $x\in N$ with $t_N(x)<T_x$, there exists $t\leqslant t_N(x)$ such that $\Phi(t)x\in E$.\eenu

\bd A pair of closed subsets $(N,\,E)$ of $X$ is called a {\bf Wa\.{z}ewski pair}, if
\benu\item[(1)]$E$ is an exit set of $N$; and
\item[(2)] $\ol{N\sm E}$ is strongly admissible.\eenu\ed

Given a subset $A$ of $X$, we always denote by $\cI(A)$ the maximal invariant set in $A$. When $A$ is strongly admissible, $\cI(A)$ is compact (see \cite{R80}).
Let $(N,E)$ be a Wa\.zewski pair with $H=\ol{N\sm E}$. We say $(N,E)$ is a {\em Wa\.zewski pair of $\cI(H)$} if it is necessary to emphasize the compact invariant set $\cI(H)$.

\bd\label{de3.2} Let $(N,E)$ be a closed pair and $H=\ol{N\sm E}$.
We say that $(N,E)$ is {\bf transversal} for $\Phi$, if for every $x\in H\cap E$ and $t>0$, there is $s\in(0,t)$ such that $\Phi(s)x\notin H\cap E$.\ed
\br Conley index pairs (\cite{Con78,R80}) and shape index pairs (\cite{WLD15,WLD20}) are common examples of Wa\.zewski pairs.
Typical transversal Wa\.zewski pairs include each isolating block with its exit set (see \cite{Con78,R80}) and even every closed invariant set with the empty set.\er
There are natural simple results about transversal Wa\.zewski pairs as follows.
\bl\label{l2.11} Suppose that a Wa\.zewski pair $(N,E)$ is transversal.
Then
\benu\item[(1)]$\cI(H)\cap E=\es$; and\item[(2)]$t_H$ is continuous on $N\sm W^s_N(\cI(H))$.\eenu
\el
\bo The first conclusion is a simple deduction.
For the second one, given $x\in N\sm W^s_N(K)$, we easily know that $t_H(x)<\8$.
Then the continuity is given by Wa\.zewski's theorem (see, \cite{Con78}).
\eo

Given a Wa\.zewski pair $(N,E)$, we now consider the pointed space $(N/E,[E])$.
The {\em quotient space} $N/E$ is defined as follows.
If $E\neq\es$, then the space $N/E$ is obtained by collapsing $E$ to a single point $[E]$ in $N\cup E$.
If $E=\es$, we choose a single isolated point $*\notin N$ and define $N/E$ to be the space $N\cup\{*\}$ equipped with the sum topology.
In the latter case we still use the notation  $[E]$ to denote the base point $*$.

Define a {\em quotient flow} $\~\Phi$ of $\Phi$ on $N/E$ as follows:

If $\~{x}=[E]$, then $$\~{\Phi}(t)\~x\equiv\~x$$ for  $t\in\R^+$; and if $\~{x}=[x]$ for some $x\in\ol{N\setminus E}$, then
$$\~{\Phi}(t)\~{x}=\left\{\ba{ll}\,[\Phi(t)x],&\mb{ for }t<t_{\ol{N\setminus E}}(x);\\[1ex]\,[E],&\mb{ for }t\geqslant t_{\ol{N\setminus E}}(x).\ea\right.$$
Since $E$ is $N$-positively invariant, it can be easily seen that $\~{\Phi}(t)$ is a well defined semigroup on $N/E$. Moreover, $\~\Phi$ is a continuous on $\R^+\X N/E$ and $N/E$ is strongly admissible for $\~\Phi$.

\bl\label{le3.5} Let $(N,E)$ be a Wa\.zewski pair and $K=\cI(\ol{N\sm E})$.
Let $\~\Phi$ be the quotient flow on $N/E$. Define for $\~\Phi$ on $N/E$
$$\cA:\,=W^u([K])\cup\{[E]\}.$$
Then $\cA$ is the compact global attractor for $\~\Phi$ in $N/E$ satisfying
\be\label{3.1}\cA=[W^u_N(K)]\cup\{[E]\}.\ee
Moreover, if $K\cap E=\es$, the attractor $\cA$ has a Morse decomposition $\sM:\,=\{\{[E]\},[K]\}$.
\el
\bo The first conclusion is similar to the results of \cite[Lemmas 4.6 and 4.7]{WLD15}.
The proof is just a slight modification of the proofs of the two lemmas, by using the framework of topological dynamical systems (\cite{LWX17}).
We thus omit it.

Now we only prove that for the case when $K\cap E=\es$, the attractor has the Morse decomposition $\sM$.
In this case, it follows from \cite{WLD15} that $[E]$ is an attractor.
We only need to show $[K]=\cA\sm\W([E])$.

Indeed, on one hand, it is obvious that $[K]\ss\cA\sm \W([E])$.
On the other hand, for every $\~x\in\cA\sm\W([E])$, there is a full solution $\~\gam$ not containing $[E]$ in $\cA$ such that $\~\gam(0)=\~x$ (otherwise $\~x\in\W([E])$).
By the definition of quotient flow, there is a full solution $\gam$ of $\Phi$ in $N\sm E$ such that $[\gam(t)]=\~\gam(t)$ for all $t\in\R$.
So $\gam$ is contained in $K$, which means $\~\gam$ is contained in $[K]$.
Hence $\~x=\~\gam(0)\in[K]$. The proof is complete.\eo
\subsection{Lyapunov functions of quotient flows}

We first introduce some typical functions that will be used in the following discussions.

Let $K\ss X$ be a closed subset and $U$ be a subset of $X$ with $K\ss U$.
A nonnegative function $\zeta\in \cC(U)$ is called a {\em $\cK_0$ function} of $K$ on $U$, if
$$\zeta(x)=0\Longleftrightarrow x\in K.$$ If moreover the {\em level set}
$$\zeta^a=\{x\in\W:\,\zeta(x)\leqslant a\}$$
is closed in $X$ for every $a\geqslant0$, we say $\zeta$ is a {\em $\cK_0^\8$ function} of $\cA$ on $\W$.

If $X$ is a metric space and $A$ is a nonempty closed subset of $X$, then the distance $\di(x,A)$ is a $\cK_0^\8$ function of $A$ on $X$.
If $B$ is another nonempty closed subset of $X$ with $A\cap B=\es$, then the function defined as
$$\frac{\di(x,A)}{\min\{1,\di(x,B)\}},\hs x\in X\sm B$$
is a $\cK_0^\8$ function of $A$ on $X\sm B$.
Thus we conclude a simple lemma.

\bl\label{le3.6} Let $A$ be a closed subset and $U$ be an open subset of a metric space $X$ with $A\ss U$.
Then, there is a $\cK_0^\8$ function $\zeta$ of $A$ on $U$ such that $\zeta(x)\geqslant\di(x,A)$ for each $x\in U$.\el

Let $\cA$ be an attractor and $\W:\,=\W(\cA)$ be the region of attraction of $\cA$.
A nonnegative continuous function $\zeta:\,\W\ra\R^+$ is said to be a {\em Lyapunov function} of $\cA$, if $\zeta$ is a $\cK_0$ function of $\cA$ on $\W$, and $\zeta(\Phi(t)x)< \zeta(x)$ for each $x\in\W\sm\cA$ and $t>0$.

Let $(N,E)$ be a Wa\.zewski pair and $\~\Phi$ be the quotient flow on $N/E$.
Then according to our previous paper \cite{WLD20}, we know that every attractor $\cA$ of $\~\Phi$ has a $\cK_0^\8$ Lyapunov function on $\W(\cA)$.
Actually, we have much more information for this Lyapunov function.

\bl\label{le3.7} Every attractor $\cA$ of $\~\Phi$ has a $\cK_0^\8$ Lyapunov function $\zeta$ on $\W(\cA)$ such that, for each $a\geqslant0$,
\be\label{3.2}\pi^{-1}(\zeta^a)\ss\mB(A,a),\ee
where $A=\pi^{-1}(\cA)$, $\pi:N\cup E\ra N/E$ is the quotient map and $\mB(A,a)$ is the set of all points in $X$ with the distance from $A$ less that $a$.\el
\bo The construction of $\zeta$ can be referred to the proof of Theorem 2.6 in \cite{WLD20}.
What we need to do is to check \eqref{3.2}.
For this we are necessary to recall some necessary constructions of $\zeta$.

If $\ol{N\sm E}\cap E=\es$ and $\cA=\{[E]\}$, we have that $\W(\cA)=\{[E]\}$.
The function $\psi([E])=0$ is just what we desire and satisfies \eqref{3.2}.

Now we only consider the case when $N\cap E\ne\es$ or $\cA\ne[E]$.

Let $U=\pi^{-1}(\W(\cA))$.
By Lemma \ref{le3.6}, we have a $\cK_0^\8$ function $\de$ of $A$ on $U$ such that $\de(x)\geqslant\di(x,A)$.
Define a function $\psi:\,\W(\cA)\ra\R^+$ such that $\psi([E])=0$ and
$$\psi(\~x)=\de(x)\mb{ for }\~x=\pi(x)\in\W(\cA)\mb{ with }x\in U\sm E.$$
It is obvious that for each $a\geqslant0$,
\be\label{3.3}\pi^{-1}(\psi^a)\ss\de^a\ss\mB(A,a).\ee
For every $\~x\in\W(\cA)$, define
$$\xi(\~x)=\sup_{t\geqslant0}\psi(\~\Phi(t)\~x)\hs\mb{and}\hs\zeta(\~x)=\xi(\~x)+\int_0^{\8}e^{-t}\xi(\~\Phi(t)\~x)dt.$$
Then $\zeta$ is the $\cK_0^\8$ Lyapunov function required (see \cite[Theorem 3.4]{WLD20}).
For every $\~x\in\W(\cA)$,
\be\label{3.4}\psi(\~x)\leqslant\xi(\~x)\leqslant\zeta(\~x).\ee
Combining (\ref{3.3}) with (\ref{3.4}), we can easily obtain \eqref{3.2}.\eo

\section{Relative Category and Morse Decomposition}\label{s3}
\subsection{Relative category}
In this section we recall the concept of relative category (see \cite{FLR94,W96}).
Let $X$ be a topological space and $I=[0,1]$.

A closed subset $A$ is {\em contractible} in $X$, if there exists $h\in \cC(I\X A,X)$, the set of all continuous maps from $I\X A$ to $X$, such that, for every $u$, $v\in A$,
$$h(0,u)=u,\hs h(1,u)=h(1,v).$$
\bd\label{de4.1} Let $A$, $B$, $Y$ be closed subsets of $X$. Then by definition, $A\prec_YB$ in $X$ if $Y\ss A\cap B$ and there exists $h\in\cC(I\X A,X)$ such that
\benu\item[(1)] $h(0,x)=x$, $h(1,x)\in B$, for all $x\in A$, and \item[(2)] $h(s,Y)\ss Y$, for all $s\in I$.\eenu\ed

\bd Let $Y\ss A$ be closed subsets of $X$. The {\bf category} of $A$ in $X$ relative to $Y$ is the least $n\in\N^+\cup\{\8\}$ such that there exists $n+1$ closed subsets $A_0$, $A_1$, $\cdots$, $A_n$ of $X$ satisfying
\benu\item[(1)]$\disp A=\Cup_{j=0}^{n}A_j$,\item[(2)]$A_1$, $\cdots$, $A_n$ are contractible in $X$, and\item[(3)] $A_0\prec_YY$ in $X$.\eenu

We denote the category of $A$ in $X$ relative to $Y$ by $\cat_{X,Y}(A)$. The {\bf category} of $A$ in $X$ is defined by $\cat_X(A):=\cat_{X,\es}(A)$.\ed

Let $A$, $B$, $Y$ be closed subsets of $X$ such that $Y\ss A$. The relative category has the following basic properties:
\benu\item[(1)] Normalisation: $\cat_{X,Y}(Y)=0$,\item[(2)] Subadditivity: $\cat_{X,Y}(A\cup B)\leqslant\cat_{X,Y}(A)+\cat_X(B)$,\item[(3)] Homotopy: if $A\prec_YB$ then $\cat_{X,Y}(A)\leqslant\cat_{X,Y}(B)$,
\item[(4)] Monotonicity: if $A\ss B\ss X$, then
$$\cat_{X,Y}(A)\leqslant\cat_{B,Y}(A)\hs\mb{and}\hs \cat_{X,Y}(A)\leqslant\cat_{X,Y}(B).$$\eenu

A metric space $X$ is an {\em absolute neighborhood extensor}, shortly an {\em ANE}, if, for every metric space $E$, every closed subset $F$ of $E$ and every map $f:F\ra X$, there exists a continuous extension of $f$ defined on a neighborhood of $F$ in $E$.
Important examples of ANE are {\em closed convex subsets of normed spaces, Banach manifolds, manifolds with boundary and finite product of ANEs} (\cite{FLR94}).

\bp[{\cite[Proposition 2.9]{FLR94}}]\label{p4.5} Let $Y$ be a closed subset of $X$ and suppose that both $X$, $Y$ are ANEs.
Then for an arbitrary closed subset $A\ss X$, there exists a closed neighborhood $B$ of $A$ such that $$\cat_{X,Y}(B)=\cat_{X,Y}(A).$$
\ep
The property given by Proposition \ref{p4.5} is called the {\em continuity property}.

\subsection{Relative category and Morse decomposition}

Let $X$ be a complete metric space with the metric $\di(\cdot,\cdot)$ and $\Phi$ be a local semiflow on $X$.
Let $K$ be a compact invariant set of $\Phi$ in $X$ and $\sM=\{M_1,\cdots,M_n\}$ be its Morse decomposition.
We have the following relation of L-S categories between $K$ and $\sM$.

\bt\label{th4.6} Let $N$ be an ANE such that $K\subset N$.
Then
\be\label{4.1}\cat_{N}(K)\leqslant\sum_{i=1}^n\cat_N(M_i).\ee
\et

\bo To show \eqref{4.1}, we first consider the case when $n=2$.
In this case $M_1$ is an attractor in $K$.
By Proposition \ref{p4.5}, we find closed neighborhood $B_i$ of $M_i$ in $N$ such that
\be\label{4.2}\cat_N(B_i)=\cat_N(M_i).\ee
Then $\ol{K\sm M_2}$ is contained in the region of attraction of $M_1$ and therefore, there exists $T>0$ such that
$$B'_1:=\Phi(T)\ol{K\sm B_2}\subset B_1\hs\mb{and}\hs\Phi(t)\ol{K\sm B_2}\subset K,\hs\mb{for }t\in[0,T],$$
which means that $\ol{K\sm B_2}\prec_{\es} B'_1$.
Combining these information and \eqref{4.2}, we obtain
\begin{align*}\cat_N(K)\leqslant&\cat_N(B_2)+\cat_N(\ol{K\sm B_2})\leqslant\cat_N(M_2)+\cat_N(B'_1)\\
\leqslant&\cat_N(M_2)+\cat_N(B_1)=\cat_N(M_1)+\cat_N(M_2),\end{align*}
which is \eqref{4.1} in case when $n=2$.

In the case when $n\geqslant3$, recalling the definition of Morse decomposition (Definition \ref{de2.4}), we have a sequence of attractors $A_i$ such that $\{A_i,M_{i+1}\}$ is a Morse decomposition of $A_{i+1}$, $i=1,\cdots,n-1$.
Hence
\begin{align*}\cat_N(K)\leqslant&\cat_N(M_n)+\cat_N(A_{n-1})\leqslant\cat_N(M_n)+\cat_N(M_{n-1})+\cat_N(A_{n-2})\\
\leqslant&\cdots\leqslant\sum_{i=1}^n\cat_N(M_i),\end{align*}
which is \eqref{4.1} in case when $n\geqslant 3$.
The proof of this theorem is complete.
\eo
Theorem \ref{th4.6} is a concise and general relation of the L-S categories of a compact invariant set and its Morse decomposition.
However, the topological structure of compact invariant sets is usually complicated, which brings great difficulty to the calculation of the corresponding categories.
Therefore, we are motivated to study more about this relation.

Now let $(N,E)$ be a transversal Wa\.zewski pair of $K$ for $\Phi$ and $H=\ol{N\sm E}$.
In the following discussion of this section, we always assume $E\ss N$.
We are to study the relation of relative categories of the pair $(N,E)$ and the Morse decomposition $\sM$.
We are hence devoted to show the following main theorem next.
\bt[Main Theorem]\label{th4.7}Suppose that $N$ is an ANE.
Then
\be\label{4.3}\cat_{N,E}(N)\leqslant\sum_{i=1}^n\cat_N(M_i).\ee\et
First, we provide some auxiliary constructions for semiflows and some necessary results.
\Vs

\bl\label{le4.8} There exists a closed neighborhood $F$ of $W^u_N(K)\cup E$ in $N$ such that $(N,F)$ and $(F,E)$ are transversal Wa\.zewski pairs of $\es$ and $K$, respectively.
Moreover,
\be\label{4.4}\cat_{N,E}(N)=\cat_{N,E}(F).\ee
\el
\bo We consider the quotient flow $\~\Phi$ on $N/E$.
By Lemma \ref{le3.5}, we know that $\{\{[E]\},[A]\}$ is a Morse decomposition of the global attractor $\cA$ of $\~\Phi$.
So $\cA$ has a $\cK_0^{\8}$ Lyapunov function $\zeta$ on $N/E$.

For each $a>0$, define $F^a=\pi^{-1}(\zeta^a)$, where $\pi:N\ra N/E$ is the quotient map.
Since $\zeta^a$ is a closed neighborhood of $\cA$, $F^a$ is a closed neighborhood of $W^u_N(K)\cup E$ in $N$.
Also since $\zeta$ is a $\cK_0^{\8}$ Lyapunov function $\zeta$ on $N/E$, then clearly $F^a$ is an exit set of $N$ and $(N,F^a)$ is a transversal Wa\.zewski pair of $\es$.
Moreover, $(F^a,E)$ is a transversal Wa\.zewski pair of $K$.
We see that each $F^a$ can be chosen to be the $F$ we desire for $a>0$.

Now we show \eqref{4.4}.
If $F=N$ (this is possible, when $N$ is contained in the unstable manifold $W^u(K)$ of $K$), \eqref{4.4} holds obviously and we only consider the case when $F\subsetneq N$.

In this case, due to that $\cI(\ol{N\sm F})=\es$, it is easy to see that $t_{\ol{N\sm F}}(x)<\8$ for every $x\in N$.
By Lemma \ref{l2.11}, $t_{H}$ is continuous on $N$.
For $s\in[0,1]$ and $x\in N$, let
$$h(s,x)=\Phi(st_{H}(x))x.$$
Then one can easily check that $h$ is a strong deformation retraction of $N$ onto $F$, which means $N\prec_{E} F$.
Hence by homotopy and monotonicity of relative category, we have
$$\cat_{N,E}(N)\leqslant\cat_{N,E}(F)\leqslant\cat_{N,E}(N),$$
which implies \eqref{4.4}.
Now the proof is finished.\eo

Let $C$ be a closed neighborhood of $K$ in $H$ with $C\cap E=\es$. Define a set
\be\label{4.5}C^s_N=\{x\in N:\,\mb{there is }t\geqslant0\mb{ such that }\Phi([0,t])x\ss N\mb{ and }\Phi(t)x\in C\}.\ee
Then we have the following consequences.
\bl\label{le4.9}\benu\item[(1)] $C^s_N\cap E=\es$ and $C^s_N$ is closed.\item[(2)] $W^s_N(K)\ss{\rm int}_{N}(C_N^s)$.\eenu\el
\bo We only consider the case when $C\ne\es$, since Lemma \ref{le4.9} clearly holds if $C=\es$.

(1) Since $C\cap E=\es$, it is clear that $C_N^s\cap E=\es$.
Let $x_n\in C^s_N$ be a sequence such that $x_n\ra x_0$.
It is sufficient to show $x_0\in C^s_N$ for the closedness.
By definition (\ref{4.5}), we have $t_n\geqslant0$ such that $\Phi(t_n)x_n\in C$.
If $t_n$ is bounded, we can assume $t_n\ra t_0$ and then by the closedness of $C$, $\Phi(t_n)x_n\ra\Phi(t_0)x_0\in C$.
This implies that $x_0\in C_N^s$.
If $t_n$ is unbounded, we can assume $t_n\ra\8$ and then by the admissibility (see \cite{R80}), $x_0\in W^s_N(K)$.
Note that $W^s_N(K)\ss C_N^s$.
Therefore $C_N^s$ is closed.
\vs

(2) We only need to show that every point $x\in W^s_N(K)$ allows an open neighborhood $U$ in $N$ such that $U\ss C_N^s$.
If this is not true, then there is a sequence $N\ni x_n\ra x$ such that for all $t\geqslant0$ satisfying $\Phi([0,t])x_n\ss N$, the end point $\Phi(t)x_n$ is not contained in $C$.

If $t_H(x_n)$ is bounded, we can assume $t_H(x_n)\ra t_0$ and then
$$\Phi(t_H(x_n))x_n\ra\Phi(t_0)x\in E,$$
which means $t_H(x)\leqslant t_0$. But it follows from $x\in W^s_N(K)$ that $t_N(x)=\8$, a contradiction.

If $t_H(x_n)$ is unbounded, we can assume $t_H(x_n)\ra\8$.
Then for every $t>0$, we have $t_H(x_n)>0$ when $n$ is large enough.
Hence we have $\Phi(t)x_n\ra\Phi(t)x\notin{\rm int}_N(C)$.
However, since $x\in W^s_N(K)$, then $w(x)\ss K$.
This indicates that there is $t>0$ such that $\Phi(t)x\in{\rm int}_N(C)$, which is also a contradiction and ends the proof.\eo

\bl\label{le4.10} Let $U$ be an open neighborhood of $K$.
Then there is a closed subset $F$ of $N$ such that $(F,E)$ is a Wa\.zewski pair of $K$, and a closed neighborhood $C$ of $K$ in $H$ with $C\cap E=\es$ such that $C^s_F\ss U$.\el

\bo If $K=\es$, let $F=N$ and $C=\es$ and we are done. Hence we assume $K\ne\es$ in the following.

Note that $K$ itself is an attractor in $K$.
We follow the proof of Lemma \ref{le4.8}.
By Lemma \ref{le3.7}, setting $A'=\pi^{-1}(\cA)$, we have $F^a\ss\mB(A',a)$, where $\cA$ is the global attractor of $\~\Phi$ in $N/E$.
By \eqref{3.1}, one has
$$\cA=[W_{F^a}^u(K)]\cup\{[E]\}\hs\mb{and}\hs A'=W_{H^a}^u(K)\cup E,$$
where $H^a=\ol{F^a\sm E}$.
By the properties of strong admissibility (see \cite{R80}), $W_{H^a}^u(K)$ is compact and $W^s_{F^a}(K)$ is closed.
By Proposition \ref{p4.5}, we can take an open neighborhood $B$ of $K$ in $U$ such that
$$2\de:=\di(B,E)>0.$$

Suppose the conclusion does not hold true.
Then we have two positive sequences $\ve_n\ra0$ and $a_n\ra0$ with $\mB(K,\ve_n)\ss B$ and $a_n<\de$ such that there is $x_n\in(\ol{\mB}_H(K,\ve_n))_{F^{a_n}}^s\sm B$, where $\ol{\mB}_H(K,\ve)=\ol{\mB}(K,\ve)\cap H$.
By the definition (\ref{4.5}), there is $t_n\geqslant0$ such that
$$\Phi([0,t_n])x_n\cap B=\es\hs\mb{and}\hs y_n:=\Phi(t_n)x_n\in\pa B.$$
This indicates that
$$\ba{rcl}y_n&\in&\pa B\cap (\ol{\mB}_H(K,\ve_n))_{F^{a_n}}^s\\&\stac{\eqref{3.2}}{\ss}&\mB(A',a_n)\sm\mB(E,\de)\\&\ss&\(\mB(W_{H^\de}^u(K),a_n)\cup\mB(E,a_n)\)\sm\mB(E,\de)\\&\ss&\mB(W_{H^\de}^u(K),a_n)\ea$$
Hence $y_n$ has a convergent subsequence (still denoted by $y_n$) such that
\be\label{4.6}y_n\ra y_0\in W_{H^\de}^u(K)\cap\pa B.\ee
Note still $y_n\in(\ol{\mB}_H(K,\ve_n))_{F^{a_n}}^s\sm\mB(K,\ve_n)$.
By (\ref{4.5}) again, there is $s_n\geqslant0$ such that
$$\Phi([0,s_n])y_n\ss F^{a_n}\hs\mb{and}\hs z_n:=\Phi(s_n)y_n\in\pa\mB(K,\ve_n),$$
and so $z_n$ can be assumed to converge to $z_0\in K$.

If $s_n$ is bounded, we can assume $s_n\ra s_0$ and then $\Phi(s_n)y_n\ra\Phi(s_0)y_0=z_0\in K$.
This means $\Phi([0,\8))y_0\ss H^\de$ and $\w(y_0)\ss A$.
So $y_0\in W^s_{H^\de}(K)$; if $s_n$ is unbounded, we can assume $s_n\ra\8$, and by the admissibility, we also have $y_0\in W^s_{H^\de}(K)$.
Recalling (\ref{4.6}), we have $y_0\in K\cap\pa B$, a contradiction!\eo

\bl\label{le4.11}Suppose that $N$ is an ANE.
Then
$$\cat_{N,E}(N)\leqslant\cat_{N}(K).$$\el

\bo By the continuity property of relative category, there is a closed neighborhood $B$ of $K$ in $N$ such that
\be\label{4.7}\cat_{N}(K)=\cat_{N}(B).\ee
By Lemma \ref{le4.10}, there is a closed subset $F$ of $N$ such that $(F,E)$ is a Wa\.zewski pair of $K$ and a closed neighborhood $C$ of $K$ in $F$ such that $C_{F}^s\subset B\sm E$.
Note $C_F^s$ is closed. Hence
\be\label{4.8}\cat_{N}(C_{F}^s)\leqslant\cat_{N}(B).\ee

Let $E'=\ol{F\sm C_{F}^s}$.
Since $W^u_N(K)\cap E'=\es$ by Lemma \ref{le4.9}, similar to the discussion in the proof of Lemma \ref{le4.8}, we have $E'\prec_EE$ in $N$.
This indicates that
\be\label{4.9}\cat_{N,E}(E')\leqslant\cat_{N,E}(E)=0.\ee
Then by the property of relative category, Lemma \ref{le4.8}, \eqref{4.7}, \eqref{4.8} and \eqref{4.9},
$$\cat_{N,E}(N)=\cat_{N,E}(F)\leqslant\cat_{N,E}(E')+\cat_{N}(C_{F}^s)\leqslant\cat_N(K),$$
which completes the proof.
\eo

Eventually, Theorem \ref{th4.7} is the direct deduction of Theorem \ref{th4.6} and Lemma \ref{le4.11}.

\bco\label{co4.12} Under the conditions of Theorem \ref{th4.7}, we have the following consequences,
\begin{align}&\cat_{N,E}(N)\leqslant\sum_{i=1}^n\cat_{M_i}(M_i)\hs\mb{and}\label{4.10}\\
&\cat_{N/E,[E]}(N/E)\leqslant\sum_{i=1}^n\cat_{M_i}(M_i).\label{4.11}\end{align}\eco
\bo For $i=1,\cdots,n$, by the monotonicity of the relative category, we know that
\be\label{4.12}\cat_N(M_i)\leqslant\cat_K(M_i)\leqslant\cat_{M_i}(M_i).\ee
The relation \eqref{4.10} follows immediately from \eqref{4.3} and \eqref{4.12}.
Thus only the relation \eqref{4.11} is not a trivial result and we prove it in the following.

Since the quotient space $N/E$ may not be an ANE, we first study the original pair $(N,E)$.
As an open subset of $N$, $N\sm E$ is an ANE.
Then by Proposition \ref{p4.5}, we have a neighborhood $B$ of $K$ in $N\sm E$ such that
\be\label{4.13}\cat_{N\sm E}(B)=\cat_{N\sm E}(K).\ee
Obviously $B$ is also a neighborhood of $K$ in $N$.
Recalling Lemmas \ref{le4.8} and \ref{le4.10}, we obtain a closed subset $G$ of $N$, such that $(N,G)$ and $(G,E)$ are transversal Wa\.zewski pairs with $\cI(G\sm E)=K$,
and a closed neighborhood $C$ of $K$ in $G$, such that $C\cap E=\es$ and $C_G^s\subset B$.
This as well as \eqref{4.13} implies
\be\label{4.14}\cat_{N\sm E}(C^s_G)\leqslant\cat_{N\sm E}(K).\ee
Moreover, by Lemma \ref{le4.9}, $C_N^s\subset N\sm E$.
By the definition \eqref{4.5}, one easily sees that $(C_N^s,C_G^s)$ is a transversal Wa\.zewski pair of $\es$, due to the fact that $(N,G)$ is a transversal Wa\.zewski pair.
By Lemma \ref{le4.8}, we know that $C_G^s$ is a strong deformation retract of $C_N^s$.
Following the homotopy property of relative category, we have
\be\label{4.15}\cat_{N\sm E}(C_N^s)\leqslant\cat_{N\sm E}(C_G^s).\ee

Now we proceed in the quotient space $N/E$ to prove it.
Let $\pi:N\ra N/E$ be the quotient map.
Since $\pi|_{N\sm E}:N\sm E\ra\pi(N\sm E)$ is a homeomorphism, it is clear that
\be\label{4.16}\cat_{\pi(N\sm E)}((\pi(C_N^s)))=\cat_{N\sm E}(C_N^s).\ee
Lemmas \ref{le3.5} and \ref{le4.9} tells us that $\pi(\ol{N\sm C_N^s})$ is contained in the region of attraction $\W([E])$ of the attractor $[E]$.
Also the singleton $[E]$ is a strong deformation retract of $\W([E])$ (see \cite[Proposition 2.5]{WLD15}).
This indicates that $\pi(\ol{N\sm C_N^s})\prec_{\{[E]\}}\{[E]\}$.
As a consequence, combining the estimates \eqref{4.14}, \eqref{4.15}, \eqref{4.16} and Theorem \ref{th4.6},
we obtain
\begin{align*}\cat_{N/E,[E]}(N/E)\leqslant&\cat_{N/E,[E]}(\pi(\ol{N\sm C_N^s}))+\cat_{N/E}(\pi(C_N^s))\\
\leqslant&\cat_{\pi(N\sm E)}(\pi(C_N^s))=\cat_{N\sm E}(C_N^s)\\
\leqslant&\cat_{N\sm E}(K)\leqslant\sum_{i=1}^n\cat_{N\sm E}(M_i)\\
\leqslant&\sum_{i=1}^n\cat_{M_i}(M_i),\end{align*}
which proves \eqref{4.11} and finishes the proof.
\eo
\br\label{re4.13} The relations \eqref{4.10} and \eqref{4.11} are a general version of Po\'zniak \cite[Theorem 3.1]{P94}.

If the Wa\.zewski pair $(N,E)$ is a Conley index pair of $K$, i.e., $\ol{N\sm E}$ is an isolating neighborhood of $K$ in $X$, a well-known result is that the the homotopy type of the quotient space $N/E$, called the Conley index of $K$, does not depend on the choice of $(N,E)$.
Then the homotopy property of relative category helps to define relative category on the Conley index of $K$.
Thus \eqref{4.11} provides us a relation between the relative categories of the Conley index and Morse decomposition of a compact isolated invariant set,
which immediately processes the continuation property under small perturbation.\er

\br Based on the development of L-S category in shape theory (see \cite{B78}), it is also possible to develop the relative category in shape theory, in which case, these consequences would have a shape-theoretical version, as another generalization of \cite{S00}.

Particularly, if we consider L-S category $\ul{\eta}$ of a metric compactum $X$ in the sense of Borsuk (\cite{B78,S00}), in the cases when $\Phi$ is a flow on $X$ or two-sided on the unstable manifold $W^u(K)$ of $K$, the following inequality (see \cite[Theorem 1 and Lemma 4]{S00}) can be obtained as well for infinite-dimensional dynamical systems,
\be\label{4.16A}\ul{\eta}(W^u(K))=\ul{\eta}(K)\leqslant\sum_{i=1}^n\ul{\eta}(M_i),\ee
with the arguments in Sanjurjo \cite{S00,S03}.\er

\section{Some Extensions of the Categories for Morse Decompositions}
In this section we discuss some topics of dynamical systems, in which the relation of (L-S, relative) category and Morse decomposition can be successfully applied.
Still, we consider a semiflow $\Phi$ on a complete metric space $X$.

First we present a consequence of detecting connecting trajectories between Morse sets in infinite dimensional dynamical systems, as a generalization of \cite[Corollary 5]{S00}.

Let $K$ be a compact invariant set with a Morse decomposition $\sM=\{M_1,\cdots,M_n\}$.
Since we do not know if $K$ is an ANE, generally it is hard to compare the (L-S) categories $\cat_K(K)$ and $\disp\sum_{i=1}^n\cat_{M_i}(M_i)$ as \eqref{4.16A}.
However, if there exist no connecting trajectories between any two Morse sets in $\sM$, we can easily obtain that
$$\cat_K(K)=\sum_{i=1}^n\cat_{M_i}(M_i).$$
As a converse-negative sentence of this result, we have the following one to detect the existence of connecting trajectories between Morse sets for semiflows on infinite-dimensional spaces.
\bt\label{th5.1} If the L-S categories of $K$ and $M_i$ satisfy the inequality
$$\cat_K(K)\neq\sum_{i=1}^n\cat_{M_i}(M_i),$$
then there exists a connecting trajectory $\gam$ between $M_i$ and $M_j$ in $K$, for some $i,j\in[1,n]$ with $i<j$, such that
$$\w(\gam)\subset M_i\hs\mb{and}\hs\al(\gam)\subset M_j.$$\et
\Vs

Now suppose that the local semiflow $\Phi$ is gradient, i.e., there is a continuous function $V:X\ra\R$, usually also called {\em Lyapunov function}, such that $t\ra V(\Phi(t)x)$ is non-increasing for each $x\in X$, and if $x$ is such that $V(\Phi(t)x)=V(x)$ for all $t\geqslant0$, then $x$ is an equilibrium of $\Phi$.
Each compact invariant set of $\Phi$ contains at least one equilibrium.
Consequently by using Theorem \ref{th4.7}, we infer the following conclusion.

\bt\label{co4.14}
Let $\Phi$ be a gradient semiflow on $X$ and $(N,E)$ be a transversal Wa\.zewski pair for $\Phi$.
Suppose that $N$ is an ANE.
Then $N\sm E$ contains at least $\cat_{N,E}(N)$ (or $\cat_{N/E,[E]}(N/E)$) equilibria of $\Phi$.\et
\bo Let $\sE(A)$ denote all the equilibria of $\Phi$ in each subset $A$ of $X$ and $\#(A)$ denote the number of points in $A$.
Let $K=\cI(\ol{N\sm E})$.

If $\#(\sE(K))=\8$, we are done.
Hence we only consider $\#(\sE(K))<\8$.
For the gradient system $\Phi$, one has $\sE(N\sm E)=\sE(K)$.
According to \cite[Proposition 5.16]{CLR13}, we impose an order on $\sE(K)$ such that $\sE(K)=\{e_1,\cdots,e_n\}$ is a Morse decomposition of $K$. Thus the conclusion follows immediately from Corollary \ref{co4.12}, since $\cat_{N}(\{e_i\})=1$ for $i=1,\cdots,n$.\eo

\br\label{re4.15} Consider the solutions of the following equation
\be\label{4.17}-\De u+f(x,u(x))=0,\hs\mb{in }U;\Hs u(x)=0\hs\mb{on }\pa U,\ee
where $U\subset\R^n$ is a bounded domain and $f$ is a nonlinear function.

The corresponding evolution equation of \eqref{4.17}
\be\label{4.18}\frac{\di u}{\di t}+Au+f(x,u)=0\ee
generates a $C_0$ gradient semiflow $\Phi$ on $H_0^1(\W)$.
Let
$$\vp(u)=\frac12\int_U\|\na u\|^2\di x+\int_U\int_0^uf(s)\di s\di x.$$
Then $\vp$ is the variational functional of \eqref{4.17} and a Lyapunov function of $\Phi$, which is thus gradient.
Note that the critical points of $\vp$ are the solutions of \eqref{4.18} and the equilibria of $\Phi$, and also that (PS)$_c$ (PS represents Palais-Smale) condition for $\vp$ is somehow equivalent to the strong admissibility condition for $\Phi$.
As a result, critical point theorem of $\vp$ with relative category (see \cite[Theorem 5.19]{W96}) coincides with Corollary \ref{co4.14} by choosing appropriate Wa\.zewski pairs.

In this sense, Corollary \ref{co4.14} can be viewed as a simple generalization of critical point theorems with (L-S, relative) category in dynamical systems.
\er
\bibliographystyle{plain}

\end{document}